\documentclass{amsart}

\usepackage{amsthm} 
\usepackage{amsmath}
\usepackage{amssymb} 

\usepackage{bussproofs}

\EnableBpAbbreviations
\newcommand{\AXM}[1]{\AXC{$#1$}}
\newcommand{\UIM}[1]{\UIC{$#1$}}
\newcommand{\BIM}[1]{\BIC{$#1$}}

\newcommand{\limplies}{\rightarrow}
\newcommand{\liff}{\leftrightarrow}
\newcommand{\ex}[1]{\exists #1 \;} 
\newcommand{\fa}[1]{\forall #1 \;} 

\newtheorem{theorem}{Theorem}[section]
\newtheorem{lemma}[theorem]{Lemma}
\newtheorem{proposition}[theorem]{Proposition}

\newcommand{\na}[1]{\mathit{#1}}    
\newcommand{\ax}[1]{\mathit{(#1)}}  

\newcommand{\lam}[1]{\lambda #1 \;} 
\newcommand{\proves}{\vdash}

\newcommand{\la}{(}
\newcommand{\ra}{)}
\newcommand{\st}{ \; | \; } 
\newcommand{\ph}{\varphi}
\newcommand{\concat}{\mathord{\hat{\;}}}
\newcommand{\dash}{\mathalpha{\mbox{-}}} 

\newcommand{\tsub}{\mathbin{\mathchoice
{\buildrel .\lower.6ex\hbox{\vphantom{.}} \over {\smash-}}%
{\buildrel .\lower.6ex\hbox{\vphantom{.}} \over {\smash-}}%
{\buildrel .\lower.4ex\hbox{\vphantom{.}} \over {\smash-}}%
{\buildrel .\lower.3ex\hbox{\vphantom{.}} \over {\smash-}}}}

\title[Functional interpretation and inductive definitions]{Functional
  interpretation and\\ inductive definitions}

\thanks{Work by the first author has been partially supported by NSF
  grant DMS-0700174 and a grant from the John Templeton Foundation.}

\author{Jeremy Avigad and Henry Towsner}

\begin{document}

\begin{abstract}
  Extending G\"odel's \emph{Dialectica} interpretation, we provide a
  functional interpretation of classical theories of positive
  arithmetic inductive definitions, reducing them to theories of
  finite-type functionals defined using transfinite recursion on
  well-founded trees.
\end{abstract}

\maketitle


\section{Introduction}
\label{introduction:section}

Let $X$ be a set, and let $\Gamma$ be a monotone operator from the
power set of $X$ to itself, so that $A \subseteq B$ implies $\Gamma(A)
\subseteq \Gamma(B)$. Then the set
\[
I = \bigcap \{ A \st \Gamma(A) \subseteq A \}
\]
is a least fixed point of $\Gamma$; that is, $\Gamma(I) = I$, and $I$
is a subset of any other set with this property. $I$ can also be
characterized as the limit of a sequence indexed by a sufficiently
long segment of the ordinals, defined by $I_0 = \emptyset$,
$I_{\alpha+1} = \Gamma(I_\alpha)$, and $I_\lambda = \bigcup_{\gamma <
  \lambda} I_\gamma$ for limit ordinals $\gamma$. Such inductive
definitions are common in mathematics; they can be used, for example,
to define substructures generated by sets of elements, the collection
of Borel subsets of the real line, or the set of well-founded trees on
the natural numbers.

From the point of view of proof theory and descriptive set theory, one
is often interested in structures that are countably based, that is,
can be coded so that $X$ is a countable set. In that case, the
sequence $I_\alpha$ stabilizes before the least uncountable ordinal.
In many interesting situations, the operator $\Gamma$ is given by a
positive arithmetic formula $\ph(x,P)$, in the sense that $\Gamma(A) =
\{ x \st \ph(x,A) \}$ and $\ph$ is an arithmetic formula in which the
predicate $P$ occurs only positively. (The positivity requirement can
be expressed by saying that no occurrence of $P$ is negated when $\ph$
is written in negation-normal form.)

The considerations above show that the least fixed point of a positive
arithmetic inductive definition can be defined by a $\Pi^1_1$ formula.
An analysis due to Stephen Kleene~\cite{kleene:55,kleene:55a} shows
that, conversely, a positive arithmetic inductive definition can be
used to define a complete $\Pi^1_1$ set. In the 1960's, Georg Kreisel
presented axiomatic theories of such inductive definitions
\cite{kreisel:63,buchholz:et:al:81}. In particular, the theory
$\na{ID_1}$ consists of first-order arithmetic augmented by additional
predicates intended to denote least fixed-points of positive
arithmetic operators.  $\na{ID_1}$ is known to have the same strength
as the subsystem $\na{\Pi^1_1\dash CA}^-$ of second order arithmetic,
which has a comprehension axiom asserting the existence of sets of
numbers defined by $\Pi^1_1$ formulas without set parameters. It also
has the same strength as Kripke Platek admissible set theory,
$\na{KP\omega}$, with an axiom asserting the existence of an infinite
set. (See \cite{buchholz:et:al:81,jaeger:86} for details.)

A $\Pi_2$ sentence is one of the form $\fa{\bar x} \ex {\bar y} R(\bar
x, \bar y)$, where $\bar x$ and $\bar y$ are tuples of variables
ranging over the natural numbers, and $R$ is a primitive recursive
relation. Here we are concerned with the project of characterizing the
$\Pi_2$ consequences of the theories $\na{ID_1}$ in constructive or
computational terms. This can be done in a number of ways. For
example, every $\Pi_2$ theorem of $\na{ID_1}$ is witnessed by a
function that can be defined in a language of higher-type
functionals allowing primitive recursion on the natural numbers as
well as a schema of recursion along well-founded trees, as described
in Section~\ref{theories:section} below. We are particularly
interested in obtaining a translation from $\na{ID_1}$ to a
constructive theory of such functions that makes it possible to ``read
off'' a description of the witnessing function from the proof of a
$\Pi_2$ sentence in $\na{ID_1}$.

There are currently two ways of obtaining this information. The first
involves using ordinal analysis to reduce $\na{ID_1}$ to a
constructive analogue \cite{buchholz:81b,pohlers:75,pohlers:81}, such
as the theory $\na{ID_1^{i,sp}}$ discussed below, and then using
either a realizability argument or a Dialectica interpretation of the
latter \cite{howard:72,buchholz:81}. One can, alternatively, use a
forcing interpretation due to Buchholz \cite{buchholz:81,aehlig:03} to
reduce $\na{ID_1}$ to $\na{ID_1^{i,sp}}$.

Here we present a new method of carrying out this first step, based on
a functional interpretation along the lines of G\"odel's
``Dialectica'' interpretation of first-order arithmetic. Such
functional interpretations have proved remarkably effective in
``unwinding'' computational and otherwise explicit information from
classical arguments (see, for example,
\cite{kohlenbach:05,kohlenbach:to:appear,kohlenbach:oliva:03a}).
Howard \cite{howard:72} has provided a functional interpretation for a
restricted version of the constructive theory $\na{ID_1^{i,sp}}$, but
the problem of obtaining such an interpretation for classical theories
of inductive definitions is more difficult, and was posed as an
outstanding problem in \cite[Section 9.8]{avigad:feferman:98}.
Feferman \cite{feferman:68} used a Dialectica interpretation to obtain
ordinal bounds on the strength of $\na{ID_1}$ (the details are
sketched in \cite[Section 9]{avigad:feferman:98}), and
Zucker~\cite{zucker:73} used a similar interpretation to bound the
ordinal strength of $\na{ID_2}$.  But these interpretations do not
yield $\Pi_2$ reductions to constructive theories, and hence do not
provide computational information; nor do the methods seem to extend
to the theories beyond $\na{ID_2}$. Our interpretation bears
similarities to those of Burr \cite{burr:00} and Ferreira and Oliva
\cite{ferreira:oliva:05}, but is not subsumed by either; some of the
differences between the various approaches are indicated in
Section~\ref{functional:interp:section}.
  
The outline of the paper is as follows. In
Section~\ref{theories:section}, we describe the relevant theories and
provide an overview of our results. Our interpretation of $\na{ID_1}$
is presented in three steps. In Section~\ref{embedding:section}, we
embed $\na{ID_1}$ in an intermediate theory, $\na{OID_1}$,
which makes the transfinite construction of the fixed-point explicit.
In Section~\ref{functional:interp:section}, we present a functional
interpretation that reduces $\na{OID_1}$ to a second
intermediate theory, $\na{Q_0 T_\Omega + \ax{I}}$. Finally, the latter
theory is interpreted in a constructive theory, $\na{QT_\Omega^i}$,
using a cut elimination argument in Section~\ref{cut:elim:section}. In
Section~\ref{iterating:section}, we show that our interpretation
extends straightforwardly to cover theories of iterated inductive
definitions as well.

We are grateful to Solomon Feferman, Philipp Gerhardy, Paulo Oliva,
and Wilfried Sieg for feedback on an earlier draft; and we are
especially grateful to Fernando Ferreira for a very careful reading
and substantive corrections.


\section{Background}
\label{theories:section}

In this paper, we interpret classical theories of inductively defined
sets in constructive theories of transfinite recursion on well-founded
trees. In this section, we describe the relevant theories, and provide
an overview of our results.

Take classical first-order Peano arithmetic, $\na{PA}$, to be
formulated in a language with symbols for each primitive recursive
function and relation. The axioms of $\na{PA}$ consist of basic axioms
defining these functions and relations, and the schema of induction,
\[
\ph(0) \land \fa x (\ph(x) \limplies \ph(x+1)) \limplies \fa x \ph(x),
\]
where $\ph$ is any formula in the language, possibly with free
variables other than $x$. $\na{ID_1}$ is an extension of $\na{PA}$
with additional predicates $I_\psi$ intended to denote the least fixed
point of the positive arithmetic operator given by
$\psi$. Specifically, let $\psi(x,P)$ be an arithmetic formula with at
most the free variable $x$, in which the predicate symbol $P$ occurs
only positively. We adopt the practice of writing $x \in I_\psi$
instead of $I_\psi(x)$. $\na{ID_1}$ then includes the following
axioms:
\begin{itemize}
\item $\fa x (\psi(x,I_\psi) \limplies x \in I_\psi)$
\item $\fa x (\psi(x,\theta/P)\limplies\theta(x)) \limplies
\fa {x \in I_\psi} \theta(x)$, for each formula $\theta(x)$.
\end{itemize}
Here, the notation $\psi(\theta/P)$ denotes the result of replacing
each atomic formula $P(t)$ with $\theta(t)$, renaming bound variables
to prevent collisions. The first axiom asserts that $I_\psi$ is closed
with respect to $\Gamma_\psi$, while the second axiom schema expresses
that $I_\psi$ is the smallest such set, among those sets that can be
defined in the language. Below we will use the fact that this schema,
as well as the schema of induction, can be expressed as rules. For
example, $I_\psi$-leastness is equivalent to the rule ``from $\fa x
(\psi(x,\theta'/P) \limplies \theta'(x))$ conclude $\fa {x \in I_\psi}
\theta'(x)$.'' To see this, note that the rule is easily justified
using the corresponding axiom; conversely, one obtains the axiom for
$\theta(x)$ by taking $\theta'(x)$ to be the formula $(\fa z
(\psi(z,\theta/P) \limplies \theta(z))) \limplies \theta(x)$ in the
rule.

One can also design theories of inductive definitions based on
intuitionistic logic. In order for these theories to be given a
reasonable constructive interpretation, however, one needs to be more
careful in specifying the positivity requirement on $\psi$. One option
is to insist that $P$ does not occur in the antecedent of any
implication, where $\lnot \eta$ is taken to abbreviate $\eta \limplies
\bot$. Such a definition is said to be \emph{strictly positive}, and
we denote the corresponding axiomatic theory $\na{ID_1^{i,sp}}$.  An
even more restrictive requirement is to insist that $\psi(x)$ is of
the form $\fa {y \prec x} P(y)$, where $\prec$ is a primitive
recursive relation. These are called {\em accessibility} inductive
definitions, and serve to pick out the well-founded part of the
relation. In the case where $\prec$ is the ``child-of'' relation on a
tree, the inductive definition picks out the well-founded part of that
tree. We will denote the corresponding theory $\na{ID_1^{i,acc}}$.

The following conservation theorem can be obtained via an ordinal
analysis \cite{buchholz:et:al:81} or the methods of
Buchholz~\cite{buchholz:81}:
\begin{theorem}
\label{idone:conserve:thm}
Every $\Pi_2$ sentence provable in $\na{ID_1}$ is provable in
$\na{ID_1^{i,acc}}$. 
\end{theorem}
The methods we introduce here provide another route to this result.

Using a primitive recursive coding of pairs and writing $x \in I_y$
for $\la x, y \ra \in I$ allows us to code any finite or infinite
sequence of sets as a single set. One can show that in any of the
theories just described, any number of inductively defined sets can
coded into a single one, and so, for expository convenience, we will
assume that each theory uses only a single inductively defined set. 

We now turn to theories of transfinite induction and recursion on
well-founded trees. The starting point is a quantifier-free theory,
$\na{T_\Omega}$, of computable functionals over the natural numbers
and the set of well-founded trees on the natural numbers. In
particular, $\na{T_\Omega}$ extends G\"odel's theory $\na{T}$ of
computable functionals over the natural numbers. We begin by reviewing
the theory $\na{T}$. The set of \emph{finite types} is defined
inductively, as follows:
\begin{itemize}
\item $N$ is a finite type; and
\item assuming $\sigma$ and $\tau$ are finite types, so are $\sigma
  \times \tau$ and $\sigma \to \tau$.
\end{itemize}
In the ``full'' set-theoretic interpretation, $N$ denotes the set of
natural numbers, $\sigma \times \tau$ denotes the set of ordered pairs
consisting of an element of $\sigma$ and an element of $\tau$, and
$\sigma \to \tau$ denotes the set of functions from $\sigma$ to
$\tau$. But we can also view the finite types as nothing more than
datatype specifications of computational objects.  The set of
\emph{primitive recursive functionals of finite type} is a set of
computable functionals obtained from the use of explicit definition,
application, pairing, and projections, and a scheme allowing the
definition of a new functional $F$ by primitive recursion:
\[
\begin{split}
F(0) & = a \\
F(x+1) & = G(x,F(x))
\end{split}
\]
Here, the range of $F$ may be any finite type. The theory $\na{T}$
includes defining equations for all the primitive recursive
functionals, and a rule providing induction for quantifier-free
formulas $\ph$:
\begin{prooftree}
\AXM{\ph(0)}
\AXM{\ph(x) \limplies \ph(S(x))}
\BIM{\ph(t)}
\end{prooftree}
G\"odel's \emph{Dialectica} interpretation shows:
\begin{theorem}
\label{pa:witness:thm}
If $\na{PA}$ proves a $\Pi_2$ theorem $\fa {\bar x} \ex {\bar y}
R(\bar x, \bar y)$, there is a sequence of function symbols $\bar f$
such that $T$ proves $R(\bar x, \bar f(\bar x))$. In particular, every
$\Pi_2$ theorem of $\na{PA}$ is witnessed by sequence of primitive
recursive functionals of type $N^k \to N$.
\end{theorem}
See \cite{goedel:58,avigad:feferman:98,troelstra:90} for details. If
$(st)$ is used to denote the result of applying $s$ to $t$, we adopt
the usual conventions of writing, for example, $stuv$ for
$(((st)u)v)$. To improve readability, however, we will also sometimes
adopt conventional function notation, and write $s(t,u,v)$ for the
same term.

In order to capture the $\Pi_2$ theorems of $\na{ID_1}$, we use an
extension of $\na{T}$ that is essentially due to Howard
\cite{howard:72}, and described in \cite[Section
9.1]{avigad:feferman:98}. Extend the finite types by adding a new base
type, $\Omega$, which is intended to denote the set of well-founded
(full) trees on $N$. We add a constant, $e$, which denote the tree
with just one node, and two new operations: $\sup$, of type $(N
\to \Omega) \to \Omega$, which forms a new tree from a sequence of
subtrees, and $\sup^{-1}$, of type $\Omega \to (N \to \Omega)$,
which returns the immediate subtrees of a nontrivial tree. We extend
the schema of primitive recursion on $N$ in $T$ to the larger
system, and add a principle of primitive recursion on $\Omega$:
\[
\begin{split}
F(e) &= a\\
F(\sup \; h) &= G(\lam n F(h(n))),
\end{split}
\]
where the range of $F$ can be any of the new types. We call the
resulting theory $\na{T_\Omega}$, and the resulting set of functionals
the \emph{primitive recursive tree functionals}. Below we will adopt
the notation $\alpha[n]$ instead of $\sup^{-1}(\alpha,n)$ to denote the
$n$th subtree of $\alpha$. In that case definition by transfinite
recursion can be expressed as follows:\footnote{We are glossing over
  issues involving the treatment of equality in our descriptions of
  both $\na{T}$ and $\na{T_\Omega}$. All of the ways of dealing with
  equality in $\na{T}$ described in \cite[Section
  2.5]{avigad:feferman:98} carry over to $\na{T_\Omega}$, and our
  interpretations work with even the most minimal version of equality
  axioms associated with the theory denoted $\na{T_0}$ there. In
  particular, our interpretations to not rely on extensionality, or
  the assumption $\fa n (\alpha[n] = \beta[n]) \limplies \alpha =
  \beta$. 

  Our theory $\na{T_\Omega}$ is essentially the theory $\na{V}$ of
  Howard \cite{howard:72}. Our theory $\na{QT_\Omega^i}$ is
  essentially a finite-type version of the theory $\na{U}$ of
  \cite{howard:72}, and contained in the theory $\na{V^*}$ described
  there. One minor difference is that Howard takes the nodes of his
  trees to be labeled, with end-nodes labeled by a positive natural
  number, and internal nodes labeled $0$.}
\[
F(\alpha) =
\left\{
  \begin{array}{ll}
    a & \mbox{if $\alpha = e$} \\
    G(\lam n F(\alpha[n])) & \mbox{otherwise.}
  \end{array}
\right.
\]
A trick due to Kreisel (see \cite{howard:68,howard:72}) allows us to
derive a quantifier-free rule of transfinite induction on $\Omega$ in
$\na{T_\Omega}$, using induction on $N$ and transfinite recursion.
\begin{proposition}
\label{kreisel:prop}
The following is a derived rule of $\na{T_\Omega}$:
\begin{prooftree}
  \AXM{\ph(e,x)} 
  \AXM{\alpha \neq e \land \ph(\alpha[g(\alpha,x)],
    h(\alpha,x)) \limplies \ph(\alpha,x)} 
  \BIM{\ph(s,t)}
\end{prooftree} 
for quantifier-free formulas $\ph$. 
\end{proposition}
For the sake of completeness, we sketch a proof in the Appendix.

We define $\na{Q T_\Omega}$ to be the extension of $\na{T_\Omega}$
which allows quantifiers over all the types of the latter theory, 
strengthening the previous transfinite induction rule with a full
transfinite induction axiom schema,
\[
\ph(e) \land \fa \alpha (\alpha \neq e \land \fa n \ph(\alpha[n])
\limplies \ph(\alpha)) \limplies \fa \alpha \ph(\alpha)
\]
where $\ph$ is any formula in the expanded language. Let
$\na{QT_\Omega^i}$ denote the version of this theory based on
intuitionistic logic.

We can also add to $\na{QT_\Omega^i}$ an \emph{$\omega$-bounding}
axiom schema,
\[
\fa x \ex \alpha \psi(x,\alpha) \limplies \ex \beta \fa x \ex i
\psi(x,\beta[i]),
\]
where $x$ is of type $N$, $\alpha$ is of type $\Omega$, and $\psi$
only has quantifiers over $N$. The following theorem shows that all of
the intuitionistic theories described in this section are ``morally
equivalent,'' and reducible to $\na{T_\Omega}$.

\begin{theorem}
\label{intuitionistic:theories:strength:thm:a}
The following theories all prove the same $\Pi_2$ sentences:
\begin{enumerate}
\item $\na{ID_1^{i,sp}}$
\item $\na{ID_1^{i,acc}}$
\item $\na{QT_\Omega^i + \ax{\omega\dash bounding}}$
\item $\na{QT_\Omega^i}$
\item $\na{T_\Omega}$
\end{enumerate}
\end{theorem}

\begin{proof}
  Buchholz~\cite{buchholz:81} presents a realizability interpretation
  of $\na{ID_1^{i,sp}}$ in the theory $\na{ID_1^{i,acc}}$.
  Howard~\cite{howard:72} presents an embedding of $\na{ID_1^{i,acc}}$
  in $\na{QT_\Omega^i + \ax{\omega\dash
      bounding}}$. Howard~\cite{howard:72} also presents a functional
  interpretation of $\na{QT_\Omega^i + \ax{\omega\dash bounding}}$ in
  $\na{T_\Omega}$, which is included in $\na{QT_\Omega^i}$;
  Proposition~\ref{kreisel:prop} is used to interpret the transfinite
  induction axioms of the source theory. Interpreting $\na{T_\Omega}$
  in $\na{ID_1^{i,sp}}$ is straightforward, using the set $O$ of
  Church-Kleene ordinal notations to interpret the type $\Omega$, and
  interpreting the constants of $\na{T_\Omega}$ as hereditarily
  recursive operations over $O$ (see \cite[Sections 4.1, 9.5, and
  9.6]{avigad:feferman:98}).
\end{proof}

In fact, Howard's work \cite{howard:72} shows that
Theorem~\ref{intuitionistic:theories:strength:thm:a} still holds as
stated if one allows arbitrary formulas $\psi(x,\alpha)$ of
$\na{QT^i_\Omega}$ in the $\omega$-bounding axiom schema. In the
classical theories considered below, however, the restriction to
arithmetic quantifiers is necessary. We have therefore chosen to use
the name $\ax{\omega\dash bounding}$ for the restricted version.

We can now describe our main results. In
Sections~\ref{embedding:section} to \ref{cut:elim:section}, we present
the interpretation outlined in the introduction, which yields:
\begin{theorem}
\label{idone:witness:thm:a}
Every $\Pi_2$ sentence provable in $\na{ID_1}$ is provable in
$\na{QT_\Omega^i}$.
\end{theorem}
In fact, if $\na{ID_1}$ proves a $\Pi_2$ theorem $\fa {\bar x} \ex
{\bar y} R(\bar x, \bar y)$, our proof yields a sequence of function
symbols $\bar f$ such that $\na{QT_\Omega^i}$ proves $R(\bar x, \bar
f(\bar x))$. By Theorem~\ref{intuitionistic:theories:strength:thm:a},
this last assertion can even be proved in $\na{T_\Omega}$. Thus we
have:
\begin{theorem}
\label{idone:witness:thm}
If $\na{ID_1}$ proves a $\Pi_2$ theorem $\fa {\bar x} \ex {\bar y}
R(\bar x, \bar y)$, there is a sequence of function symbols $\bar f$
such that $\na{T_\Omega}$ proves $R(\bar x, \bar f(\bar x))$. In
particular, every $\Pi_2$ theorem of $\na{ID_1}$ is witnessed by
sequence of primitive recursive tree functionals of type $N^k \to N$.
\end{theorem}
The reduction described by Sections~\ref{embedding:section}
to \ref{cut:elim:section} is thus analogous to the reduction of
$\na{ID_1}$ given by Buchholz \cite{buchholz:81}, but relies on a
functional interpretation instead of forcing.


\section{Embedding  $\na{ID_1}$ in $\na{OID_1}$}
\label{embedding:section}

In this section, we introduce a theory $\na{OID_1}$, which makes the
transfinite construction of the fixed points of $\na{ID_1}$ explicit.
We then show that $\na{ID_1}$ is easily interpreted in $\na{OID_1}$.
The theory $\na{OID_1}$ is closely related to Feferman's theory
$\na{OR_1^\omega}$, as described in \cite{feferman:68} and
\cite[Section 9]{avigad:feferman:98}, and the embedding is similar to
the one described there.

Fix any instance of $\na{ID_1}$ with inductively defined predicate $I$
given by the positive arithmetic formula $\psi(x,P)$. The
corresponding instance of $\na{OID_1}$ is two-sorted, with variables
$\alpha, \beta, \gamma, \ldots$ ranging over type $\Omega$, and
variables $i, j, k, n, x, \ldots$ ranging over $N$. We include symbols
for the primitive recursive functions on $N$, and a function symbol
$\sup^{-1}(\alpha,n)$ which returns an element of type $\Omega$. As
above, we write $\alpha[n]$ for $\sup^{-1}(\alpha,n)$. Recall that
$\alpha[n]$ is intended to denote the $n$th subtree of $\alpha$, or
$e$ if $\alpha = e$. The language includes an equality symbol for
terms of type $N$, but \emph{not} for terms of type $\Omega$.  We
include, however, a unary predicate ``$\alpha = e$,'' which holds when
$\alpha$ is the tree with just one node. Finally, there is a binary
predicate $I(\alpha,x)$, where $\alpha$ ranges over $\Omega$ and $x$
ranges over $N$. We will write $x \in I_\alpha$ instead of
$I(\alpha,x)$, and write $x \in I_{\prec \alpha}$ for $\ex i (x \in
I_{\alpha[i]})$. The axioms of $\na{OID_1}$ are as follows:
\begin{enumerate}
\item defining axioms for the primitive recursive functions
\item induction on $N$
\item transfinite induction on $\Omega$
\item $\alpha = e \limplies \alpha[i] = e$
\item the schema of $\omega\dash\mathit{bounding}$: 
\[
\fa x \ex \alpha \ph(x,\alpha)
  \limplies \ex \beta \fa x \ex i \ph(x,\beta[i]),
\]
where $\ph$ has no quantifiers over type $\Omega$.
\item $\fa x (x \not\in I_e)$
\item $\fa \alpha (\alpha \neq e \limplies \fa x (x \in I_\alpha \liff
  \psi(x, I_{\prec \alpha})))$
\end{enumerate}
The last two axioms assert that $I_\alpha$ is the hierarchy of sets
satisfying $I_e = \emptyset$ and $I_\alpha = \Gamma_\psi(I_{\prec
  \alpha})$ when $\alpha \neq e$. For any formula $\ph$ of
$\na{ID_1}$, let $\hat \ph$ be the formula obtained by interpreting $t
\in I$ as $\ex \alpha (t \in I_\alpha)$.

\begin{theorem}
\label{main:or:thm}
If $\na{ID_1}$ proves $\ph$, then $\na{OID_1}$ proves $\hat \ph$.
\end{theorem}

We need two lemmas. In the first, let $I_{\prec \alpha} \subseteq
I_{\prec \beta}$ abbreviate the formula $\fa x (x \in
I_{\prec \alpha} \limplies x \in I_{\prec \beta})$.

\begin{lemma}
\label{or:lemma}
$\na{OID_1}$ proves that for every $\alpha$ and $i$, $I_{\prec \alpha[i]}
\subseteq I_{\prec \alpha}$.
\end{lemma}

\begin{proof}
  Use transfinite induction on $\alpha$. If $\alpha$
  is equal to $e$, the conclusion is immediate from the fourth
  axiom. In the inductive step, suppose $x$ is in $I_{\prec
    \alpha[i]}$. Then for some $j$, $x$ is in $I_{\alpha[i][j]}$. By
  the last axiom, we have $\psi(x,I_{\prec \alpha[i][j]})$. By the
  inductive hypothesis, we have $I_{\prec \alpha[i][j]} \subseteq
  I_{\prec \alpha[i]}$, and so, by the positivity of $\psi$, we have
  $\psi(x,I_{\prec \alpha[i]})$. By the last axiom again, we have $x
  \in I_{\alpha[i]}$, and hence $x \in I_{\prec \alpha}$, as
  required. 
\end{proof}

Note that if $\eta(x,P)$ is any arithmetic
formula involving a new predicate symbol $P$ and $\theta(y)$ is any
formula, applying the $\hat \cdot$-translation to $\eta(x,\theta/P)$
changes only the instances of $\theta$. In particular,
$\widehat{\eta(x,I)}$ is $\eta(x, \ex \alpha(y \in I_\alpha) / P)$.

\begin{lemma}
\label{pos:equiv:lemma}
Let $\eta(x,P)$ be a positive arithmetic formula. Then $\na{OID_1}$
proves that $\widehat{\eta(x,I)}$ implies $\ex \alpha
\eta(x,I_{\prec \alpha})$.
\end{lemma}

\begin{proof}
  Use induction on positive arithmetic formulas, expressed in
  negation-normal form. To handle the base case where $\eta(x,I)$ is
  $x \in I$, suppose we have $x \in I_\beta$. By a trivial instance of
  $\omega$-bounding, there are an $\alpha$ and an $i$ such that $x$ is
  in $I_{\alpha[i]}$. But this means that $x$ is in $I_{\prec
    \alpha}$, as required.

  All the other cases are easy, except when the outermost connective
  is a universal quantifier. In that case, suppose $\na{OID_1}$ proves
  that $\widehat{\ph(x,y,I)}$ implies $\ex \beta \ph(x,y,I_{\prec
    \beta})$. Using $\omega$ bounding, $\widehat{\fa y \ph(x,y,I)}$
  then implies $\ex \alpha \fa y \ex i \ph(x,y,I_{\prec \alpha[i]})$.
  By Lemma~\ref{or:lemma} and positivity, $\na{OID_1}$ proves $\ex
  \alpha \fa y \ph(x,y,I_{\prec \alpha})$, as required.
\end{proof}

\begin{proof}[Proof of Theorem~\ref{main:or:thm}] 
  The defining axioms for the primitive recursive functions and
  induction axioms of $\na{ID_1}$ are again axioms of $\na{OID_1}$
  under the translation, so we only have to deal with the defining
  axioms for $I$.

  To verify the translation of the closure axiom in $\na{OID_1}$,
  suppose $\widehat{\psi(x,I)}$. By Lemma~\ref{pos:equiv:lemma}, we
  have $\ex \alpha \psi(x,I_{\prec \alpha})$, which implies $\ex \alpha
  (x \in I_\alpha)$, as required.

  This leaves only the leastness property of $I$, which can be
  expressed as a rule, ``from $\fa x (\psi(x,\theta/P) \limplies
  \theta(x))$, conclude $\fa {x \in I} \theta(x)$.'' To verify the
  translation in $\na{OID_1}$, suppose $\fa x (\psi(x,\hat
  \theta/P) \limplies \hat \theta(x))$.  It suffices to show that for
  every $\alpha$, we have $\fa {x \in I_\alpha} \hat \theta(x)$. We
  use transfinite induction on $\alpha$.  In the base case, when
  $\alpha = e$, this is immediate from the defining axiom for
  $I_e$. In the inductive step, suppose we have $\fa i \fa {x \in
    I_{\alpha[i]}} \hat\theta(x)$. This is equivalent to $\fa {x \in
    I_{\prec \alpha}} \hat \theta(x)$. Using the positivity of $P$, we
  have $\fa x (\psi(x,I_{\prec \alpha}) \limplies \psi(x,\hat \theta /
  P))$. Using the definition of $I_\alpha$, we then have $\fa {x \in
    I_\alpha} \hat \theta(x)$, as required.
\end{proof}


\section{A functional interpretation of $\na{OID_1}$}
\label{functional:interp:section}

Our next step is to interpret the theory $\na{OID_1}$ in a second
intermediate theory, $\na{Q_0 T_\Omega + \ax{I}}$. First, we describe
a fragment $\na{Q_0 T_\Omega}$ of $\na{QT_\Omega}$, which is obtained
by restricting the language of $\na{QT_\Omega}$ to allow
quantification over the natural numbers only, though we continue to
allow free variables and constants of all types. We also restrict the
language so that the only atomic formulas are equalities $s = t$
between terms of type $N$. The axioms of $\na{Q_0 T_\Omega}$ are as
follows:
\begin{enumerate}
\item any equality between terms of type $N$ that can be derived in
  $\na{T_\Omega}$
\item the schema of induction on $N$.
\item the schema of transfinite induction, given as a rule:
\begin{prooftree}
\AXM{\theta(e)}
\AXM{\alpha \neq e \land \fa n \theta(\alpha[n]) \limplies
  \theta(\alpha)}
\BIM{\theta(t)}
\end{prooftree}
for any formula $\theta$ and term $t$ of type $\Omega$.
\end{enumerate}
In the transfinite induction schema, the formula $\alpha = e$ is
should be understood as the formula $f(\alpha) = 0$, where $f$ is the
function from $\Omega$ to $N$ defined recursively by $f(e) = 0, f(\sup
\; g) = 1$. Substitution is a derived rule in $\na{Q_0 T_\Omega}$,
which is to say, if the theory proves $\ph(x)$ where $x$ is a variable
of any type, it proves $\ph(s)$ for any term $s$ of that type. One can
show this by a straightforward induction on proofs, using the fact
that any substitution instance of one of the axioms or rules of
inference above is again an axiom or rule of inference. Similarly, by
induction on formulas $\ph(x)$, one can show that if $\na{T_\Omega}$
proves $s = t$ for any terms $s$ and $t$ of the appropriate type, then
$\na{Q_0 T_\Omega}$ proves $\ph(s) \liff \ph(t)$.

The following proposition shows that in $\na{Q_0 T_\Omega}$ we can use
instances of induction in which higher-type parameters are allowed to
vary. For example, the first rule states that in order to prove
$\theta(\alpha,x)$ for arbitrary $\alpha$ and $x$, it suffices to
prove $\theta(e,x)$ for an arbitrary $x$, and then, in the induction
step, prove that $\theta(\alpha,x)$ follows from
$\theta(\alpha[n],a)$, as $n$ ranges over the natural numbers and $a$
ranges over a countable sequence of parameters depending on $n$ and
$x$.

\begin{proposition}
\label{derived:rules:prop}
The following are derived rules of $\na{Q_0T_\Omega}$:
\begin{prooftree}
  \AXM{\theta(e,x)} 
  \AXM{\alpha \neq e \land \fa i \fa j
    \theta(\alpha[i], f(\alpha,x,i,j)) \limplies \theta(\alpha, x)}
  \BIM{\theta(\alpha,x)}
\end{prooftree}
and
\begin{prooftree}
  \AXM{\psi(0,x)}
  \AXM{\fa j \psi(n,f(x,n,j)) \limplies \psi(n+1,x)}
  \BIM{\psi(n,x)}
\end{prooftree}
\end{proposition}

For a fixed instance of $\na{ID_1}$, we now define the theory $\na{Q_0
  T_\Omega + \ax{I}}$ by adding a new binary
predicate $I(\alpha,x)$, which is allowed to occur in the induction
axioms and the transfinite induction rules, and the following
axioms:
\begin{enumerate}
\item[(4)] $\fa x (x \not\in I_e)$.
\item[(5)] $\fa \alpha (\alpha \neq e \limplies \fa x (x \in I_\alpha \liff
  \psi(x, I_{\prec \alpha})))$.
\item[(6)] $s \in I_\alpha \liff t \in I_\beta$ whenever 
  $s$, $\alpha$, $t$, and $\beta$ are terms such that
  $\na{T_\Omega}$ proves $s = t$ and $\alpha = \beta$.
\end{enumerate}
Proposition~\ref{derived:rules:prop} extends to this new theory, as
does the substitution rule. Thanks to axiom (6), if $\ph(x)$ is any
formula of $\na{Q_0 T_\Omega + \ax{I}}$ and $s$ and $t$ are any terms 
such that $\na{T_\Omega}$ proves $s = t$, then $\na{Q_0 T_\Omega +
  \ax{I}}$ proves $\ph(s) \liff \ph(t)$. 

The goal of this section is to use a functional interpretation to
interpret $\na{OID_1}$ in $\na{Q_0 T_\Omega + \ax{I}}$. As in
Burr~\cite{burr:00}, we use a variant of Shoenfield's interpretation
\cite{shoenfield:01} which incorporates an idea due to Diller and Nahm
\cite{diller:nahm:74}. The Shoenfield interpretation works for
classical logic, based on the connectives $\forall$, $\lor$, and
$\lnot$. This has the virtue of cutting down on the number of axioms
and rules that need to be verified, and keeping complexity
down. Alternatively, we could have used a Diller-Nahm variant of the
ordinary G\"odel interpretation, combined with a double-negation
interpretation. The relationship between the latter approach and the
Shoenfield interpretation is now well understood (see
\cite{streicher:kohlenbach:07,avigad:unp:k}).

First, we need to introduce some notation. We will often think
of an element $\alpha \neq e$ of $\Omega$ as denoting a countable set
$\{ \alpha[i] \st i \in \mathbb N \}$ of elements of $\Omega$. Within
the language of $\na{QT_\Omega}$, we therefore define $\alpha
\sqsubseteq \beta$ by
\[
\alpha \sqsubseteq \beta \equiv \fa i \ex j (\alpha[i] = \beta[j]),
\]
expressing inclusion between the corresponding sets. Let $t(i)$ be any
term of type $\Omega$, where $i$ is of type $N$. Then we can define
the union of the sets $t(0), t(1), t(2), \ldots$ by
\[
\sqcup_i t(i)\equiv \sup_j t(j_0)[j_1],
\]
where $j_0$ and $j_1$ denote the projections of $j$ under a primitive
recursive coding of pairs. In other words, $\sqcup_i t(i)$ represents
the set $\{ t(i)[k] \st i \in \mathbb N, k \in \mathbb N \}$. In
particular, we have that for every $i$, $t(i) \sqsubseteq \sqcup_i
t(i)$, since for every $k$ we have $t(i)[k] = (\sqcup_i
t(i))[(i,k)]$. Binary unions, $s \sqcup t$, can be defined in a
similar way.

We can extend these notions to higher types. Define the set of
\emph{pure $\Omega$-types} to be the smallest set of types containing
$\Omega$ and closed under the operation taking $\sigma$ and $\tau$ to
$\sigma \to \tau$. Note that every pure $\Omega$-type $\tau$ has the
form $\sigma_1 \to \sigma_2 \to \ldots \sigma_k \to \Omega$. We can
therefore lift the notions above to $a$, $b$, and $t$ of arbitrary
pure type, by defining them to hold pointwise, as follows:
\[
\begin{split}
  a[i] & \equiv \lam x ((a x)[i]) \\
  a \sqsubseteq b & \equiv \fa i \ex j \fa x ((a x)[i] =
  (b x)[j]) \\
  \sqcup_i t & \equiv \lam x (\sqcup_i (t x)) \\
  s \sqcup t & \equiv \lam x ((s x) \sqcup (t x)),
\end{split}
\]
where in each case $x$ is a tuple of variables chosen so that the
resulting term has type $\Omega$. Thus, if $a$ is of any pure type, we
can think of $a$ as representing the countable set $\{ a[i] \st i \in
\mathbb N \}$, in which case $\sqsubseteq$ and $\sqcup$ have the
expected behavior.

Below we will be interested in the situation where $\na{T_\Omega}$ can
prove $a \sqsubseteq b$ in the sense that there is an explicit term
$f(i)$, not involving $x$, such that $\na{T_\Omega}$ proves $(a x)[i]
= (b x)[f(i)]$. Notice that when $\na{T_\Omega}$ proves $a \sqsubseteq
b$ in this sense, $\na{Q_0 T_\Omega + \ax{I}}$ proves $\ph(a[i])
\limplies \ph(b[f(i)])$ for any formula $\ph$, and hence $\ex i
\ph(a[i]) \limplies \ex j \ph(b[j])$. Notice also that $\na{T_\Omega}$
proves $t(i) \sqsubseteq \sqcup_i t(i)$, $s \sqsubseteq s \sqcup t$,
and $t \sqsubseteq s \sqcup t$ in this sense.

To each formula $\ph$ in the language of $\na{OID_1}$, we associate a
formula $\ph^S$ of the form $\fa a \ex b \ph_S(a,b)$, where $a$ and
$b$ are tuples of variables of certain pure $\Omega$-types (which are
implicit in the definitions below), and $\ph_S$ is a formula in the
language of $\na{Q_0T_\Omega + \ax{I}}$.  The interpretation is
defined, inductively, in such a way that the following monotonicity
property is preserved: whenever $\na{T_\Omega}$ proves $b \sqsubseteq
b'$ in the sense above, $\na{Q_0T_\Omega + \ax{I}}$ proves $\ph_S(a,b)
\limplies \ph_S(a,b')$. In the base case, we define
\begin{align*}
I(\alpha,t)^S & \equiv I(\alpha,t) \\
(s = t)^S & \equiv s = t
\end{align*}
In the inductive step, suppose $\ph^S$ is $\fa a \ex b \ph_S(a,b)$
and $\psi^S$ is $\fa c \ex d \psi_S(c,d)$. Then we define
\begin{align*}
  (\ph \lor \psi)^S & \equiv \fa {a,c} \ex {b,d}
  (\ph_S(a,b) \lor \psi_S(c,d)) \\
  (\fa x \ph)^S & \equiv \fa a \ex b (\fa x \ph_S(a,b)) \\
  (\fa \alpha \ph)^S & \equiv \fa {\alpha,a} \ex b \ph_S(a,b) \\
  (\lnot \ph)^S & \equiv \fa B \ex a (\ex i \lnot \ph_S(a[i],B(a[i]))).
\end{align*}
Verifying the monotonicity claim above is straightforward; the inner
existential quantifier in the clause for negation takes care of the
only case that would otherwise have given us trouble. Note in
particular the clause for universal quantification over the natural
numbers. Our functional interpretation is concerned with bounds;
because we can compute ``countable unions'' using the operator
$\sqcup$, we can view quantification over the natural numbers as
``small'' and insist that the bound provided by $b$ is independent of
$x$. Note also that if $\ph$ is a purely arithmetic formula, $\ph^S$
is just $\ph$.

The rest of this section is devoted to proving the following:
\begin{theorem}
\label{functional:interp:thm}
Suppose $\na{OID_1}$ proves $\ph$, and $\ph^S$ is the formula
$\fa a \ex b \ph_S(a,b)$. Then there are terms $b$ of
$\na{T_\Omega}$ involving at most the variables $a$ and the free
variables of $\ph$ of type $\Omega$ such that $\na{Q_0T_\Omega +
  \ax{I}}$ proves $\ph_S(a,b)$.
\end{theorem}
Importantly, the terms $b$ in the statement of the theorem do not
depend on the free variables of $\ph$ of type $N$.

As usual, the proof is by induction on derivations. The details are
similar to those in Burr~\cite{burr:00}. As in Shoenfield
\cite{shoenfield:01}, we can take the logical axioms and rules to be
the following:
\begin{enumerate}
\item excluded middle: $\lnot \ph \lor \ph$
\item substitution: $\fa x \ph(x) \limplies \ph(t)$, and $\fa \alpha
  \ph(\alpha) \limplies \ph(t)$
\item expansion: from $\ph$ conclude $\ph \lor \psi$
\item contraction: from $\ph \lor \ph$ conclude $\ph$
\item cut: from $\ph \lor \psi$ and $\lnot \ph \lor \theta$, conclude
  $\psi \lor \theta$.
\item $\forall$-introduction: from $\ph \lor \psi$ conclude $\fa x \ph
  \lor \psi$, assuming $x$ is not free in $\psi$; and similarly for
  variables of type $\Omega$
\item equality axioms
\end{enumerate}

The translation of excluded middle is
\[
\fa {B,a'} \ex {a,b'} (\ex i \lnot \ph_S(a[i],B(a[i])) \lor
\ph_S(a',b')).
\]
Given $B$ and $a'$, let $a = \sup_i a'$, so that $a[i] = a'$ for every
$i$; in other words, $a$ represents the singleton set $\{ a' \}$. Let
$b' = B(a')$. Then the matrix of the formula holds with $i = 0$.

The translation of substitution for the natural numbers is equivalent
to
\[
\fa {B,a'} \ex{a,b'} (\fa{i,x} \ph_S(x,a[i],B(a[i])) \limplies
\ph_S(t,a',b')).
\]
(In this context, ``equivalent to'' means that $\na{Q_0 T_\Omega +
  \ax{I}}$ proves that the $\cdot_S$ part of the translation is
equivalent to the expression in parentheses.) Once again, given $B$
and $a'$, letting $a = \sup_i a'$ and $b' = B(a')$ works.

Handling substitution for $\Omega$ and expansion is
straightforward. Consider the contraction rule. By the inductive hypothesis
we have terms $b = b(a,c)$ and $d = d(a,c)$ satisfying
\begin{equation*}
\ph_S(a,b(a,c)) \lor \ph_S(c,d(a,c)).
\end{equation*}
Define $f(e)$ to be $b(e,e) \sqcup d(e,e)$. Then $\na{T_\Omega}$
proves $b(e,e) \sqsubseteq f(e)$ and $d(e,e) \sqsubseteq f(e)$. By
substitution and monotonicity we have $\ph_S(e,f(e)) \lor
\ph_S(e,f(e))$, and hence $\ph_S(e,f(e))$, as required.

Consider cut. By the inductive hypothesis we have terms $b = b(a,c)$
and $d = d(a,c)$ satisfying
\begin{equation}
\label{cut:eq:a}
\ph_S(a,b(a,c)) \lor \psi_S(c,d(a,c)),
\end{equation}
and terms $a' = a'(B,e)$ and $f = f(B,e)$ satisfying
\begin{equation}
\label{cut:eq:b}
\ex i \lnot \ph_S(a'(B,e)[i],B(a'(B,e))[i]) \lor \theta_S(e,f(B,e)).
\end{equation}
We need terms $d' = d'(c',e')$ and $f' = f'(c',e')$ satisfying
\[
\psi_S(c',d'(c',e')) \lor \theta_S(e',f'(c',e')).
\]
Given $c'$ and $e'$, and the terms $b(a,c)$, $d(a,c)$, $a'(B,e)$, and
$f(B,e)$, define $B' = \lam a b(a,c')$, define $a'' = \sup_i
a'(B',e')$, and then define $d' = d(a'',c')$ and $f' = f(B',e')$.
Since $\na{T_\Omega}$ proves $B'(a'') = b(a'',c')$, from
\eqref{cut:eq:a} we have
\[
\ph_S(a'',B'(a'')) \lor \psi_S(c',d').
\]
Since $a''[i] = a'(B',e')$ for every $i$, from \eqref{cut:eq:b} we
have
\[
\lnot \ph_S(a'',B'(a'')) \lor \theta_S(e',f').
\]
Applying cut in $\na{Q_0 T_\Omega + \ax{I}}$, we have $\psi_S(c',d')
\lor \theta_S(e',f')$, as required.

The treatment of $\forall$-introduction over $N$ and $\Omega$ is
straightforward. We can take the equality axioms to be reflexivity,
symmetry, transitivity, and congruence with respect to the basic
function and relation symbols in the language. These, as well as the
defining equations for primitive recursive function symbols in the
language and the defining axioms for $I$, are verified by the fact
that for formulas whose quantifiers ranging only over $N$, $\ph^S =
\ph$.

Thus we only have to deal with the other axioms of $\na{OID_1}$, namely, $\omega$ bounding, induction on $N$, and
transfinite induction on $\Omega$. Note that if $\ph$ has quantifiers
ranging only over $N$, the definition of $\exists$ in terms of
$\forall$ implies that $(\ex \alpha \ph(\alpha))^S$ is equivalent to
$\ex \alpha \ex i \ph(\alpha[i])$. To interpret the translation of
$\omega$-bounding, we therefore need to define a term $\beta =
\beta(\alpha)$ satisfying
\[
\fa x \ex i \ph_S(x,\alpha[i]) \limplies \ex j \fa x
\ex k \ph_S(x,(\beta[j])[k]).
\]
Setting $\beta = \sup_j \alpha$ means that for every $j$ we have
$\beta[j] = \alpha$, so this $\beta$ works.

We can take induction on the natural numbers to be given by the rule
``from $\ph(0)$ and $\ph(x) \limplies \ph(x+1)$ conclude $\ph(t)$ for
any term $t$.'' From a proof of the first hypothesis, we obtain a term
$b = b(a)$ satisfying
\begin{equation}
\label{ind:eqa}
\ph_S(0,a,b).
\end{equation}
From a proof of the second hypothesis, we obtain terms $a' =
a'(B',a'')$ and $b'' = b''(B',a'')$ satisfying
\begin{equation}
\label{ind:eqb}
\fa i \ph_S(x,a'[i],B'(a'[i])) \limplies \ph_S(x+1,a'',b'').
\end{equation}
If suffices to define a function $f(x,\hat a)$ and show that we can
prove
\begin{equation}
\label{ind:eqc}
\ph_S(x,\hat a, f(x,\hat a)),
\end{equation}
since if we then define $\hat b(\hat a) = \sqcup_x f(x,\hat a)$, we
have $\ph_S(x,\hat a, \hat b)$ by the monotonicity property of our
translation. Define $f$ by
\[
\begin{split}
f(0,\hat a) &= b(\hat a) \\
f(x+1,\hat a) &= b''(\lam a f(x,a), \hat a)
\end{split}
\]
Let $B'$ denote $\lam a f(x,a)$, so $f(x+1,\hat a) = b''(B',\hat a)$.
Let $A(x,\hat a)$ denote the formula \eqref{ind:eqc}. From
\eqref{ind:eqa}, we have $A(0,\hat a)$, and from \eqref{ind:eqb} we
have $\fa i A(x, a'(B',\hat a)[i]) \limplies A(x+1,\hat a)$. Using
Proposition~\ref{derived:rules:prop}, we obtain $A(x,\hat a)$, as
required.

Transfinite induction, expressed as the rule ``from $\ph(e)$ and $\fa
n \ph(\alpha[n]) \limplies \ph(\alpha)$ conclude $\ph(\alpha)$,'' is
handled in a similar way. From a proof of the first hypothesis we
obtain a term $b = b(a)$ satisfying
\begin{equation}
\label{trans:a}
\ph_S(e,a,b).
\end{equation}
From a proof of the second hypothesis we obtain terms $a' =
a'(\alpha,B',a'')$ and $b'' = b''(\alpha,B',a'')$ satisfying
\begin{equation}
\label{trans:b}
\fa i \fa n \ph_S(\alpha[n],a'[i],B'(a'[i]))
\limplies \ph_S(\alpha,a'',b'').
\end{equation}
It suffices to define a function $f$ satisfying
\[
\ph_S(\alpha,\hat a, f(\alpha, \hat a))
\]
for every $\alpha$ and $\hat a$, since then $\hat b = f(\alpha, \hat
a)$ is the desired term. Let $A(\alpha,\hat a)$ be this
last formula, and define $f$ by recursion on $\alpha$:
\[
f(\alpha, \hat a) = \left\{
\begin{array}{ll}
b(a) & \mbox{if $\alpha = e$} \\
b''(\alpha,\lam a (\sqcup_j f(\alpha[j],a)), \hat a) &
\mbox{otherwise.}
\end{array}\right.
\]
Write $B'$ for the expression $\lam a (\sqcup_j f(\alpha[j],a))$, so
we have $f(\alpha, \hat a) = b''(\alpha,B',\hat a)$ when $\alpha \neq
e$. We will use the transfinite induction rule given by
Proposition~\ref{derived:rules:prop} to show that $A(\alpha,\hat a)$
holds for every $\alpha$ and $\hat a$. From \eqref{trans:a}, we have
$A(e,\hat a)$, so it suffices to show
\[
\alpha \neq e \land \fa {n,i} 
A(\alpha[n], a'[i]) \limplies
A(\alpha,\hat a),
\]
where $a'$ is the term $a'(\alpha,B',\hat a)$. Arguing in
$\na{Q_0T_\Omega + \ax{I}}$, assume $\alpha \neq e$ and $\fa {n, i} 
A(\alpha[n], a'[i])$, that is,
\[
\fa {n, i} \ph_S(\alpha[n],a'[i],f(\alpha[n],a'[i])).
\]
By monotonicity, we have
\[
\fa {n, i} \ph_S(\alpha[n],a'[i],
\sqcup_j f(\alpha[j],a'[i])).
\]
By the definition of $B'$, this is just
\[
\fa {n, i} \ph_S(\alpha[n], a'[i], B'(a'[i])).
\]
By \eqref{trans:b}, this implies
\[
\ph_S(\alpha, \hat a, f(\alpha,\hat a)),
\]
which is $A(\alpha, \hat a)$ as required. This concludes the proof of
Theorem~\ref{functional:interp:thm}.

Our theory $\na{Q_0 T_\Omega + \ax{I}}$ is inspired by Feferman
\cite{feferman:68}, and, in particular, the theory denoted $\na{T_\Omega +
  \ax{\mu}}$ in \cite[Section 9]{avigad:feferman:98}. That theory,
like $\na{Q_0 T_\Omega + \ax{I}}$, combines a classical treatment of
quantification over the natural numbers with a constructive treatment
of the finite types over $\Omega$. 

The principal novelty of our interpretation, however, is the use of
the Diller-Nahm method in the clause for negation, and the resulting
monotonicity property. This played a crucial rule in the
interpretation of transfinite induction. The usual Dialectica
interpretation would require us to choose a single candidate for the
failure of an inductive hypothesis, something that cannot be done
constructively. Instead, using the Diller-Nahm trick, we recursively
``collect up'' a countable sequence of possible counterexamples. (The
original Diller-Nahm trick involved using only finite sequences of
counterexamples; we are grateful to Paulo Oliva for pointing out to us 
that the extension of the method to more general sequences of
counterexamples seems to have been first used by Stein~\cite{stein:78}.)

Similar uses of monotonicity can be found in functional
interpretations developed by
Kohlenbach~\cite{kohlenbach:to:appear,kohlenbach:oliva:03a} and
Ferreira and Oliva~\cite{ferreira:oliva:05}, as well as in the forcing
interpretations described in Avigad~\cite{avigad:00b}. The functional
interpretations of Avigad~\cite{avigad:98}, Burr~\cite{burr:00}, and
Ferreira and Oliva~\cite{ferreira:oliva:05} also make use of the
Diller-Nahm trick. But Kohlenbach, Ferreira, and Oliva rely on
majorizability relations, which cannot be represented in $\na{Q_0
  T_\Omega}$, due to the restricted uses of quantification in that
theory. Our interpretation is perhaps closest to the one found in
Burr~\cite{burr:00}, but a key difference is in our interpretation of
universal quantification over the natural numbers; as noted above,
because we are computing bounds and our functionals are closed under
countable sequences, the universal quantifier is absorbed by the
witnessing functional.


\section{Interpreting $\na{Q_0 T_\Omega + \ax{I}}$ in
  $\na{QT_\Omega^i}$}
\label{cut:elim:section}

The hard part of the interpretation is now behind us. It is by now
well known that one can embed infinitary proof systems for classical
logic in the various constructive theories listed in
Theorem~\ref{intuitionistic:theories:strength:thm:a}. This idea was
used by Tait \cite{tait:70}, to provide a constructive consistency
proof for the subsystem $\na{\Sigma^1_1\dash CA}$ of second-order
arithmetic. It was later used by Sieg~\cite{sieg:77,sieg:81} to
provide a direct reduction of the classical theory $\na{ID_1}$ to the
constructive theory $\na{ID^{i,sp}_2}$, as well as the corresponding
reductions for theories of transfinitely iterated inductive
definitions (see Section~\ref{iterating:section}). Here we show that,
in particular, one can define an infinitary proof system in
$\na{QT^i_\Omega}$, and use it to interpret $\na{Q_0 T_\Omega +
  \ax{I}}$ in a way that preserves $\Pi_2$ formulas. The methods are
essentially those of Sieg~\cite{sieg:77,sieg:81}, adapted to the
theories at hand. In fact, our interpretation yields particular
witnessing functions in $\na{T_\Omega}$, yielding
Theorem~\ref{idone:witness:thm:a}.

Let us define the set of infinitary \emph{constant} propositional
formulas, inductively, as follows:
\begin{itemize}
\item $\top$ and $\bot$ are formulas.
\item If $\ph_0, \ph_1, \ph_2, \ldots$ are formulas, so are
  $\bigvee_{i \in N} \ph_i$ and $\bigwedge_{i \in N} \ph_i$.
\end{itemize}
Take a \emph{sequent} $\Gamma$ to be a finite set of such formulas. As
usual, we write $\Gamma, \Delta$ for $\Gamma \cup \Delta$ and $\Gamma,
\ph$ for $\Gamma \cup \{ \ph \}$. We define a cut-free infinitary
proof system for such formulas with the following rules:
\begin{itemize}
\item $\Gamma, \top$ is an axiom for each sequent $\Gamma$.
\item From $\Gamma, \ph_i$ for some $i$ conclude $\Gamma, \bigvee_{i
    \in N} \ph_i$.
\item From $\Gamma, \ph_i$ for every $i$ conclude $\Gamma, \bigwedge_{i
    \in N} \ph_i$.
\end{itemize}
We also define a mapping $\ph \mapsto \lnot \ph$ recursively, as
follows:
\begin{itemize}
\item $\lnot \top = \bot$
\item $\lnot \bot = \top$
\item $\lnot \bigvee_{i \in N} \ph_i = \bigwedge_{i \in N} \lnot \ph_i$.
\item $\lnot \bigwedge_{i \in N} \ph_i = \bigvee_{i \in N} \lnot \ph_i$.
\end{itemize}
Note that the proof system does not include the cut rule, namely,
``from $\Gamma, \ph$ and $\Gamma, \lnot \ph$ include $\Gamma$.'' In
this section we will show that it is possible to represent
propositional formulas and infinitary proofs in the language of
$\na{QT_\Omega^i}$ in such a way that $\na{QT_\Omega^i}$ proves that
the set of provable sequents is closed under cut. We will then show
that this infinitary proof system makes it possible to interpret
$\na{Q_0T_\Omega + \ax{I}}$ in a way that preserves $\Pi_2$
sentences. This will yield Theorem~\ref{idone:witness:thm:a}. In
fact, our interpretation will yield explicit functions witnessing the
truth of the $\Pi_2$ from the proof in $\na{Q_0T_\Omega + \ax{I}}$.

We can represent formulas in $\na{QT_\Omega^i}$ as well-founded trees
whose end nodes are labeled either $\top$ or $\bot$ and whose
internal nodes are labeled either $\bigvee$ or $\bigwedge$. A
well-founded tree is simply an element of $\Omega$. As in the Appendix,
if $\alpha$ is an element of $\Omega$, then one can assign to each
node of $\alpha$ a unique ``address,'' $\sigma$, where $\sigma$ is a
finite sequence of natural numbers. Since these can be coded as
natural numbers, a labeling of $\alpha$ from the set $\{ \top, \bot,
\bigvee, \bigwedge \}$ is a function $l$ from $N$ to $N$. The
assertion that $\alpha, l$ is a formula, i.e.~that the labeling has the
requisite properties, is given by a universal formula in
$\na{QT_\Omega^i}$. Using $\lambda$-abstraction we can define
functions $F$ with recursion of the following form:
\[
F(\alpha, l) = \left\{
\begin{array}{ll}
G(l(\emptyset)) & \mbox{if $\alpha = e$} \\
H(\lam n F(\alpha[n], \lam \sigma l(\la n \ra\concat \sigma))) &
\mbox{otherwise,}
\end{array}
\right.
\]
where $\emptyset$ denotes the sequence of length $0$. This yields a
principle of recursive definition on formulas, which can be used, for
example, to define the map $\ph \to \lnot \ph$. (This particular
function can be defined more simply by just switching $\top$ with
$\bot$ and $\bigwedge$ with $\bigvee$ in the labeling.) A principle of
induction on formulas is obtained in a similar way. We can now
represent proofs as well-founded trees labeled by finite sets of
formulas and rules of inference, yielding principles of induction and
recursion on proofs as well.

We will write $\proves \Gamma$ for the assertion that $\Gamma$ has an
infinitary proof, and we will write $\proves \ph$ instead of $\proves
\{ \ph \}$. The proofs of the following in $\na{QT_\Omega^i}$ are now
standard and straightforward (see, for example,
\cite{schwichtenberg:77alt,sieg:81}).

\begin{lemma}[Weakening]
If $\proves \Gamma$ and $\Gamma' \supseteq \Gamma$ then $\proves
\Gamma'$. 
\end{lemma}

\begin{lemma}[Excluded middle]
For every formula $\ph$, $\proves \{ \ph, \lnot \ph \}$.
\end{lemma}

\begin{lemma}[Inversion]
\mbox{}
\begin{itemize}
\item If $\proves \Gamma, \bot$, then $\proves \Gamma$.
\item If $\proves \Gamma, \bigwedge_{i \in N} \ph_i$, then $\proves
  \Gamma, \ph_i$ for every $i$.
\end{itemize}
\end{lemma}

The first and third of these is proved using induction on proofs in
$\na{QT_\Omega^i}$. The second is proved using induction on formulas.

\begin{lemma}[Admissibility of cut]
If $\proves \Gamma, \ph$ and $\proves \Gamma, \lnot \ph$, then
$\proves \Gamma$.
\end{lemma}

\begin{proof}
We show how to cast the usual proof as a proof by induction on
formulas, with a secondary induction on proofs. For any formula $\ph$,
define
\[
\ph^{\lor} = 
\left\{
\begin{array}{ll}
  \ph & \mbox{if $\ph$ is $\top$ or of the 
    form $\bigvee_{i \in N} \psi_i$} \\
  \lnot \ph & \mbox{otherwise.}
\end{array}\right.
\]
We express the claim
to be proved as follows:
\begin{quote}
For every formula $\ph$, for every proof $d$, the following holds: if
$d$ is a proof of a sequent of the form $\Gamma, \ph^{\lor}$, then
$\proves \Gamma, \lnot (\ph^{\lor})$ implies $\proves \Gamma$.
\end{quote}
The most interesting case occurs when $\ph = \ph^{\lor}$ is of the
form $\bigvee_{i \in N} \psi_i$, and the last inference of $d$ is
of the form
\begin{prooftree}
\AXM{\Gamma, \bigvee_{i \in N} \psi_i, \psi_j}
\UIM{\Gamma, \bigvee_{i \in N} \psi_i}
\end{prooftree}
Given a proof of $\Gamma, \bigwedge_{i \in N} \lnot \psi_i$, apply
weakening and the inner inductive hypothesis for the immediate
subproof of $d$ to obtain a proof of
$\Gamma, \psi_j$, apply inversion to obtain a proof of $\Gamma, \lnot
\psi_j$, and then apply the outer inductive hypothesis to the subformula
$\lnot \psi_j$ of $\ph$. 
\end{proof}

We now assign, to each formula $\ph(\bar x)$ in the language of
$\na{Q_0T_\Omega^i + \ax{I}}$, an infinitary formula $\hat \ph(\bar
x)$. More precisely, to each formula $\ph(\bar x)$ we assign a
function $F_\ph(\bar x)$ of $\na{T_\Omega}$, in such a way that
$\na{QT_\Omega^i}$ proves ``for every $\bar x$, $F_\ph(\bar x)$ is an
infinitary propositional formula.'' We may as well take $\lor$,
$\lnot$, and $\forall$ to be the logical connectives of $\na{Q_0
  T_\Omega^i + \ax{I}}$, and use the Shoenfield axiomatization of
predicate logic given in the last section. For formulas not involving
$I_\alpha$, the assignment is defined inductively as follows:
\begin{itemize}
\item $\widehat{s = t}$ is equal to $\top$ if $s = t$, and $\bot$
  otherwise.
\item $\widehat{\ph \lor \psi}$ is equal to $\bigvee_j
  \widehat{\theta_j}$, where $\theta_0 = \ph$ and $\theta_j = \psi$
  for $j > 0$. 
\item $\widehat{\fa x \ph(x)}$ is $\bigwedge_j \widehat{\ph}(j)$.
\item $\widehat{\lnot \ph}$ is $\lnot \widehat{\ph}$.
\end{itemize}
If $I$ corresponds to the inductive definition $\psi(x,P)$, the
interpretation of $x \in I_\alpha$ is defined recursively:
\[
x \in I_\alpha = \left\{
\begin{array}{ll}
\bot & \mbox{if $\alpha = e$} \\
\widehat{\psi(x,I_{\prec \alpha})} & \mbox{otherwise.}
\end{array}\right.
\]

The following lemma asserts that this interpretation is sound.

\begin{lemma}
  If $\na{Q_0T_\Omega + \ax{I}}$ proves $\ph(\bar x)$, then
  $\na{QT_\Omega^i}$ proves that for every $\bar x$, $\proves \hat
  \ph(\bar x)$.
\end{lemma}

\begin{proof}
  We simply run through the axioms and rules of inference in $\na{Q_0
    T_\Omega + \ax{I}}$. If $s = t$ is a theorem of $\na{T_\Omega}$,
  it is also a theorem of $\na{Q T_\Omega^i}$.  Hence
  $\na{QT_\Omega^i}$ proves $\widehat{s = t} = \top$, and so $\proves
  \widehat{s = t}$.

The interpretation of the logical axioms and rules are easily
validated in the infinitary propositional calculus augmented with the
cut rule, and the interpretation of the defining axioms for $I_\alpha$
are trivially verified given the translation of $\widehat{t \in
  I_\alpha}$. This leaves only induction on $N$ and transfinite
induction on $\Omega$. We will consider transfinite induction on
$\Omega$; the treatment of induction on $N$ is similar. 

We take transfinite induction to be given by the rule ``from $\ph(e)$
and $\alpha \neq e \land \fa n \ph(\alpha[n]) \limplies \ph(\alpha)$
conclude $\ph(\alpha)$.'' Arguing in $\na{Q_0T_\Omega + \ax{I}}$,
suppose for every instantiation of $\alpha$ and the parameters of
$\ph$ there is an infinitary derivation of the $\widehat{\cdot}$
translation of these hypothesis. Use
transfinite induction to show that for every $\alpha$ there is an
infinitary proof of $\widehat \ph(\alpha)$. When $\alpha = e$, this is
immediate. In the inductive step we have infinitary proofs of
$\widehat{\ph}(\alpha[n])$ for every $n$. Applying the
$\bigwedge$-rule, we obtain an infinitary proof of $\widehat{\fa n
  \ph(\alpha[n])}$, and hence, using ordinary logical operations in
the calculus with cut, a proof of $\widehat{\ph}(\alpha)$.
\end{proof}

We note that with a little more care, one can obtain cut-free proofs
of the induction and transfinite induction axioms; see, for example,
\cite{buchholz:81}.

\begin{lemma}
Let $\ph$ be a formula of the form $\fa {\bar x} \ex {\bar y} R(\bar
x,\bar y)$, where $R$ is primitive recursive. Then $\na{QT_\Omega^i}$
proves that $\proves \hat \ph$ implies $\ph$.
\end{lemma}

\begin{proof}
  Using a primitive recursive coding of tuples we can assume, without
  loss of generality, that each of $\bar x$ and $\bar y$ is a single
  variable. Using the inversion lemma, it suffices to prove the
  statement for $\Sigma_1$ formulas, which we can take to be of the
  form $\ex y S(y)$ for some primitive recursive $S$. Use induction on
  proofs to prove the slightly more general claim that given any proof
  of either $\{ \widehat{\ex y S(y)} \}$ or $\{ \widehat{\ex y S(y)},
  \bot \}$ there is a $j$ satisfying $S(j)$.  In a proof of either
  sequent, the last rule rule can only have been a $\bigvee$ rule,
  applied to a sequent of the form $\{ \widehat{S(j)} \}$ or $\{
  \widehat{\ex y S(y)}, \widehat{S(j)} \}$. If $\widehat{S(j)}$ equals
  $\top$, $j$ is the desired witness; otherwise, apply the inductive
  hypothesis.
\end{proof}

Putting the pieces together, we have shown:
\begin{theorem}
\label{last:step:thm}
Every $\Pi_2$ theorem of $\na{Q_0 T_\Omega + \ax{I}}$ is a
theorem of $\na{QT_\Omega^i}$.
\end{theorem}
Together with Theorems~\ref{main:or:thm} and
\ref{functional:interp:thm}, this yields
Theorem~\ref{idone:witness:thm:a}.  Note that every time we used
induction on formulas or proofs in the lemmas above, the arguments
give explicit constructions that are represented by terms of
$\na{T_\Omega}$. So we actually obtain, from an $\na{ID_1}$ proof of a
$\Pi_2$ sentence, a $\na{T_\Omega}$ term witnessing the conclusion and
a proof that this is the case in $\na{QT_\Omega^i}$. By
Theorem~\ref{intuitionistic:theories:strength:thm:a}, this can be
converting to a proof in $\na{T_\Omega}$, if desired.

Our reduction of $\na{ID_1}$ to a constructive theory has been carried
out in three steps, amounting, essentially, to a functional
interpretation on top of a straightforward cut elimination argument. A
similar setup is implicit in the interpretation of $\na{ID_1}$ due to
Buchholz~\cite{buchholz:81}, where a forcing interpretation is used in
conjunction with an infinitary calculus akin to the one we have used
here. We have also considered alternative reductions of $\na{Q_0
  T_\Omega + \ax{I}}$ that involve either a transfinite version of the
Friedman $A$-translation \cite{friedman:78} or a transfinite version
of the Dialectica interpretation. These yield interpretations of
$\na{Q_0 T_\Omega + \ax{I}}$ not in $\na{QT_\Omega^i}$, however, but
in a Martin-L\"of type theory $\na{ML_1V}$ with a universe and a type
of well-founded sets \cite{aczel:78}. $\na{ML_1V}$ is known to have
the same strength as $\na{ID_1}$, but although many consider
$\na{ML_1V}$ to be a legitimate constructive theory in its own right,
we do not know of any reduction of $\na{ML_1V}$ to one of the other
constructive theories listed in
Theorem~\ref{intuitionistic:theories:strength:thm:a} that does not
subsume a reduction of $\na{ID_1}$. Thus the methods described in
this section seem to provide an easier route to a stronger result.


\section{Iterating the interpretation}
\label{iterating:section}

In this section, we consider theories $\na{ID_n}$ of finitely iterated
inductive definitions. These are defined in the expected way:
$\na{ID_{n+1}}$ bears the same relationship to $\na{ID_n}$ that
$\na{ID_1}$ bears to $\na{PA}$. In other words, in $\na{ID_{n+1}}$ one
can introduce a inductive definitions given by formulas $\psi(x,P)$,
where $\psi$ is a formula in the language of $\na{ID_n}$ together with
the new predicate $P$, in which $P$ occurs only positively.

Writing $\Omega_0$ for $N$ and $\Omega_1$ for $\Omega$, we can now
define a sequence of theories $\na{T_{\Omega_n}}$. For each $n \geq 1$
take $\na{T_{\Omega_{n+1}}}$ to add to $\na{T_{\Omega_n}}$ a type
$\Omega_{n+1}$ of trees branching over $\Omega_n$, with corresponding
constant $e$ and functionals
$\sup:(\Omega_n\rightarrow\Omega_{n+1})\rightarrow\Omega_{n+1}$ and
$\sup^{-1}:\Omega_{n+1}\rightarrow(\Omega_n\rightarrow\Omega_{n+1})$.
Once again, we extend primitive recursion in $\na{T_{\Omega_n}}$ to
the larger system and add a principle of primitive recursion on
$\Omega_{n+1}$. The theories $\na{QT^i_{\Omega_{n+1}}}$ are defined
analogously.  It is convenient to act as though for each $i < j$,
$\Omega_j$ is closed under unions indexed by $\Omega_i$; this can
arranged by fixing injections of each $\Omega_i$ into $\Omega_{j-1}$.

In this section, we show that our interpretation extends to
$\na{ID_n}$, to yield the following generalization of
Theorem~\ref{idone:witness:thm:a}:
\begin{theorem}
\label{idone:witness:thm:z}
Every $\Pi_2$ sentence provable in $\na{ID_n}$ is provable in
$\na{QT_{\Omega_n}^i}$.
\end{theorem}
As with Theorem~\ref{idone:witness:thm:a}, the proof yields a
particular term witnessing the $\Pi_2$ assertion, and the correctness
of that witnessing term can be established in $\na{T_{\Omega_n}}$, by
a generalization of Theorem
\ref{intuitionistic:theories:strength:thm:a}. The interpretation can
be further extended to theories of transfinitely iterated inductive
definitions, as described in \cite{buchholz:et:al:81}. We do not,
however, know of any ordinary mathematical arguments that are
naturally represented in such theories.

To extend the theories $\na{OID_1}$ to theories $\na{OID_n}$, we first
have to generalize the schema of $\omega$-bounding. For each $i < j$,
define the schema of $\Omega_i\dash\Omega_j$-bounding as follows:
\begin{equation*}
  \fa {\alpha^{\Omega_i}} \ex {\beta^{\Omega_j}} \ph(\alpha,\beta)
  \limplies  \ex {\beta^{\Omega_j}} \fa {\alpha^{\Omega_i}} 
  \ex {\gamma^{\Omega_i}} 
  \ph(\alpha,\beta[\gamma]).
\end{equation*}
for every formula $\ph$ with quantifiers ranging over the types
$\Omega_0, \ldots, \Omega_i$. With this notation, the
$\omega$-bounding schema is now corresponds to
$\Omega_0\dash\Omega_1$-bounding.

We extend the theories $\na{OID_1}$ to theories $\na{OID_n}$ in
the expected way, where now $\na{OID_n}$ includes the schema of
$\Omega_i\dash\Omega_j$-bounding for each $i < j \leq n$. The fixed
points $I_1,\ldots,I_n$ of $\na{ID_n}$ are interpreted iteratively
according to the recipe in Section~\ref{embedding:section}. In
particular, if $\psi_j(x,P)$ is gives the definition of the $j$th
inductively defined predicate $I_j$, the translation of $\psi_j$ has
quantifiers ranging over at most $\Omega_{j-1}$, and $t \in I_j$ is
interpreted as $\ex {\alpha^{\Omega_j}} (t \in I_{j,\alpha})$, where
the predicates $I_j(\alpha,x)$ are defined in analogy to
$I(\alpha,x)$. This yields:
\begin{theorem}
If $\na{ID_n}$ proves $\ph$, then $\na{OID_n}$ proves $\hat
\ph$.
\end{theorem}

Next, we define theories $\na{Q_{n-1} T_{\Omega_n} + \ax{I}}$ in
analogy to the theory $\na{Q_0 T_\Omega + \ax{I}}$ of
Section~\ref{functional:interp:section}, except that we include the
$\Omega_i\dash\Omega_j$-bounding axioms for $i < j < n$ in
$\na{Q_{n-1} T_{\Omega_n} + \ax{I}}$. Now it is quantification over
the types $\Omega_0, \ldots, \Omega_{n-1}$ that is considered
``small,'' and absorbed into the target theory. In particular, for $i
< j < n$, the $\Omega_i\dash\Omega_j$-bounding axioms of $\na{OID_n}$
are unchanged by the functional interpretation. The
$\Omega_i\dash\Omega_n$ bounding axioms for $i < n$, induction on $N$,
and transfinite induction on $\Omega_n$ are interpreted as in
Section~\ref{functional:interp:section}. With the corresponding
modifications to $\ph^S$, we then have the analogue to
Theorem~\ref{functional:interp:thm}:
\begin{theorem}
  Suppose $\na{OID_n}$ proves $\ph$, and $\ph^S$ is the formula $\fa a
  \ex b \ph_S(a,b)$. Then there are terms $b$ of $\na{T_{\Omega_n}}$
  involving at most the variables $a$ and the free variables of $\ph$
  of type $\Omega_n$ such that $\na{Q_{n-1}T_{\Omega_n} + \ax{I}}$
  proves $\ph_S(a,b)$.
\end{theorem}

In the last step, we have to embed $\na{Q_{n-1} T_{\Omega_n} + \ax{I}}$
into an infinitary proof system in $\na{QT_{\Omega_n}^i}$. The method
of doing this is once again found in \cite{sieg:77,sieg:81}, and an
extension of the argument described in Section~\ref{cut:elim:section}.
We extend the definition of the infinitary propositional formulas so
that when, for each $\alpha\in\Omega_j$ with $j<n$, $\ph_\alpha$ is a
formula, so are $\bigvee_{\alpha\in\Omega_j}\ph_\alpha$ and
$\bigwedge_{\alpha\in\Omega_j}\ph_\alpha$. The proof of cut
elimination, and the verification of transfinite induction and the
defining axioms for the predicates $I_j(\alpha,x)$, are essentially
unchanged. The only additional work that is required is to handle the
bounding axioms; this is taken care of using a style of bounding
argument that is fundamental to the ordinal analysis of such
infinitary systems (see \cite{pohlers:75,pohlers:81,sieg:77,sieg:81}).

\begin{lemma}
  For every $i < j < n$, $\na{QT_{\Omega_n}^i}$ proves the translation
  of the $\Omega_i$-$\Omega_j$ bounding axioms.
\end{lemma}

\begin{proof}[Proof (sketch).]
  Since $\na{QT_{\Omega_n}^i}$ establishes the provability of the law
  of the excluded middle in the infinitary language, it suffices to
  show that for every sequent $\Gamma$ with quantifiers ranging over
  at most $\Omega_i$, if $\proves \Gamma, \fa {\alpha^{\Omega_i}} \ex
  {\beta^{\Omega_j}} \ph(\alpha,\beta)$, then there is a $\beta$ in
  $\Omega^j$ such that for every $\alpha$ in $\Omega_i$, $\proves \ex
  {\gamma^{\Omega_i}} \ph(\alpha,\beta[\gamma])$. But this is
  essentially a consequence of the ``Boundedness lemma for $\Sigma$''
  in Sieg~\cite[page 16 2]{sieg:81}; the requisite $\beta$ is defined
  by an explicit recursion on the derivation.
\end{proof}

This gives us the proper analogue of Theorem~\ref{last:step:thm}, and
hence Theorem~\ref{idone:witness:thm:z}.

\begin{theorem}
  Every $\Pi_2$ theorem of $\na{Q_{n-1}T_{\Omega_n} + \ax{I}}$ is a
  theorem of $\na{QT_{\Omega_n}^i}$.
\end{theorem}


\section*{Appendix: Kreisel's trick and induction with parameters}
\label{kreisel:section}

For completeness, we sketch a proof of
Proposition~\ref{kreisel:prop}. Full details can be found in
\cite{howard:68,howard:72}. 

\newtheorem*{propositiona}{Proposition~\ref{kreisel:prop}}

\begin{propositiona}
The following is a derived rule of $\na{T_\Omega}$:
\begin{prooftree}
  \AXM{\ph(e,x)} 
  \AXM{\alpha \neq e \land \ph(\alpha[g(\alpha,x)],
    h(\alpha,x)) \limplies \ph(\alpha,x)} 
  \BIM{\ph(s,t)}
\end{prooftree} 
for quantifier-free formulas $\ph$. 
\end{propositiona}

\begin{proof}
  We associate to each node of an element of $\Omega$ a finite
  sequence $\sigma$ of natural numbers, where the $i$th child of the
  node corresponding to $\sigma$ is assigned $\sigma \concat \la i
  \ra$. Then the subtree $\alpha_\sigma$ of $\alpha$ rooted at
  $\sigma$ (or $e$ if $\sigma$ is not a node of $\alpha$) can be
  defined by recursion on $\Omega$ as follows:
\[
\begin{split}
e_{\sigma} & = e \\
(\sup f)_\sigma & =
\left\{
  \begin{array}{ll}
    \sup f & \mbox{if $\sigma = \emptyset$} \\
    (f(i))_{\tau} & \mbox{if $\sigma = \la i \ra \concat \tau$}
  \end{array}
\right.
\end{split}
\]
Here $\emptyset$ denotes the sequence of length $0$. 

Now, given $\ph$, $g$, and $h$ as in the statement of the lemma, we
define a function $k(\alpha,x,n)$ by primitive recursion on $n$. The
function $k$ uses the the second clause of the rule to compute a
sequence of pairs $\la \sigma, y \ra$ with the property that
$\ph(\alpha_\sigma,y)$ implies $\ph(\alpha,x)$. For readability, we
fix $\alpha$ and $x$ and write $k(n)$ instead of $k(\alpha,x,n)$. We
also write $k_0(n)$ for $\la k(n) \ra_0$ and $k_1(n)$ for $\la k(n)
\ra_1$.
\[
\begin{split}
  k(0) &= \la \emptyset, x \ra \\
  k(n+1) &= \left\{
  \begin{array}{ll}
    \la k_0(n) \concat \la g(\alpha_{k_0(n)}, k_1(n)) \ra, \\
          \quad h(\alpha_{k_0(n)}, k_1(n)) \ra &
    \mbox{if $\alpha_{k_0(n)} \neq e$} \\
    k(n) & \mbox{otherwise.}
  \end{array}
\right.
\end{split}
\]
Ordinary induction on the natural numbers shows that for every $n$,
$\ph(\alpha_{k_0(n)},k_1(n))$ implies $\ph(\alpha,x)$.  So, it
suffices to show that for some $n$, $\alpha_{k_0(n)} = e$.

Since $k_0(0) \subseteq k_0(1) \subseteq k_0(2) \subseteq \ldots$ is
an increasing sequence of sequences, it suffices to establish the more
general claim that for every $\alpha$ and every function $f$ from $N$
to $N$, there is an $n$ such that $\alpha_{\la f(0), \ldots, f(n-1)
  \ra} = e$. To that end, by recursion on $\Omega$, define
\[
\begin{split}
g(\alpha,f) = \left\{
  \begin{array}{ll}
     1 + g(\alpha[f(0)], \lam n f(n+1) ) & \mbox{if $\alpha \neq e$}
     \\
     0 & \mbox{otherwise}
  \end{array}
\right.
\end{split}
\]
Let $h(m) = g(\alpha_{\la f(0), \ldots, f(m-1) \ra}, \lam n
f(n+m))$. By induction on $m$ we have $h(0) = m + h(m)$ as long as
$\alpha_{\la f(0),\ldots,f(m-1) \ra} \neq e$. In particular, setting
$m = h(0)$, we have $h(h(0)) = 0$, which implies $\alpha_{\la
  f(0),\ldots,f(h(0)-1) \ra} = e$, as required.
\end{proof}

The following principles of induction and recursion were used in
Section~\ref{functional:interp:section}.

\newtheorem*{propositionb}{Proposition~\ref{derived:rules:prop}}

\begin{propositionb}
The following are derived rules of $\na{Q_0T_\Omega}$:
\begin{prooftree}
  \AXM{\theta(e,x)} 
  \AXM{\alpha \neq e \land \fa i \fa j
    \theta(\alpha[i], f(\alpha,x,i,j)) \limplies \theta(\alpha, x)}
  \BIM{\theta(\alpha,x)}
\end{prooftree}
and
\begin{prooftree}
  \AXM{\psi(0,x)}
  \AXM{\fa j \psi(n,f(x,n,j)) \limplies \psi(n+1,x)}
  \BIM{\psi(n,x)}
\end{prooftree}
\end{propositionb}

\begin{proof}
  Consider the first rule. For any element $\alpha$ of $\Omega$ and
  finite sequence of natural numbers $\sigma$ (coded as a natural
  number), once again we let $\alpha_\sigma$ denote the subtree of
  $\alpha$ rooted at $\sigma$.  Let $\tau$ be the type of $x$. We will
  define a function $h(\alpha,g,\sigma)$ by recursion on $\alpha$,
  which returns a function of type $N \to \tau$, with the property
  that $h(\alpha,g,\emptyset) = g$, and for every $\sigma$,
  $\theta(\alpha_\sigma,x)$ holds for every $x$ in the range of
  $h(\alpha,g,\sigma)$. Applying the conclusion to $h(\alpha,\lam i x,
  \emptyset)$ will yield the desired result.

  The function $h$ is defined as follows:
\[
h(\alpha,g,\sigma) = 
\left\{
\begin{array}{ll}
g  & \mbox{if $\alpha = e$ or $\sigma = \emptyset$} \\
h(\alpha[i], \lam l f(\alpha,g(l_0),i,l_1), \sigma') & \mbox{if $\alpha \neq
  e$ and $\sigma = \sigma' \concat \la i \ra$}
\end{array}
\right.
\]
Using transfinite induction on $\alpha$, we have
\[
\fa \sigma \fa v \theta(\alpha_\sigma, h(\alpha,g,\sigma)(v))
\]
for every $\alpha$, and hence and hence $\theta(\alpha,h(\alpha,\lam i
x, \emptyset)(0))$. Since $h(\alpha,\lam i x, \emptyset)(0) = (\lam i
x)(0) = x$, we have the desired conclusion.

The second rule is handled in a similar way.
\end{proof}


\end{document}